\documentclass[12pt, twoside, a4paper]{article}

\usepackage{amsthm}
\usepackage{amssymb}
\usepackage{amsmath}
\usepackage{mathtools}
\usepackage{stmaryrd}

\usepackage[english]{babel}
\usepackage[utf8]{inputenc}

\usepackage[T1]{fontenc}

\usepackage{enumitem}
\setlist{nosep} 

\usepackage{tikz-cd}

\usepackage{xcolor,colortbl}

\usepackage{adjustbox}

\usepackage{caption}

\usepackage{titletoc}
\dottedcontents{section}[1.8em]{}{1.7em}{1pc}

\usepackage{hyperref}
\hypersetup{pdfborder=0 0 0}

\usepackage[all]{xy}

\usepackage{datetime}
\newdateformat{monthyeardate}{
  \monthname[\THEMONTH] \THEYEAR}

\newtheorem{theorem}{Theorem}[section]

\newtheorem{remark}[theorem]{Remark}

\newcommand{\C}{\mathbb{C}}

\newcommand{\R}{\mathbb{R}}

\newlength{\hnor}
\newlength{\dnor}
\newlength{\currenttextsize}
\makeatletter
\newcommand{\fixhd}[1]{
  \setlength{\currenttextsize}{\f@size pt}
  \raisebox{0pt}[0pt][0pt]{\fontsize{\currenttextsize}{1cm}\selectfont $#1$}
}
\makeatother

\renewcommand{\epsilon}{\varepsilon}
\usepackage{mathrsfs}

\usepackage{color,soul}

\begin{document}
\title{Maslov class of exact Lagrangians and cylindrical handles}
\author{Axel Husin}
\date{\monthyeardate\today}
\maketitle

\begin{abstract}
  A fundamental and deep result in symplectic topology due to Abouzaid and Kragh states that the Maslov class vanishes for closed exact Lagrangians in cotangent bundles of closed manifolds. In this article we prove by an explicit construction that the Maslov class does not vanish in general for closed exact Lagrangians in Weinstein domains obtained by performing a critical handle attachment to a cotangent bundle. We also define cylindrical handles as a generalization of Weinstein handles.
\end{abstract}

\tableofcontents

\section{Introduction}
In this article, we study the Maslov class of exact Lagrangians in a type of symplectic manifolds called Weinstein domains. These concepts are introduced below. More details can be found in \cite{mcduff2017introduction,geiges2008introduction,cieliebak2012stein,MR2441780,smoczyk2019generalized}.

  An exact symplectic manifold $(W,\omega=d\lambda)$ is a $2n$-dimensional manifold equipped with a one-form $\lambda\in \Omega_{\mathrm{dr}}^1W$ (called the Liouville form) such that $\omega=d\lambda$ (called the symplectic form) is non-degenerate. The vector field $Z$ defined by $\omega(Z,-)=\lambda$ is called the Liouville vector field. An $n$-dimensional submanifold $L\subset W$ such that $\lambda|_{L}\in \Omega_{\mathrm{dr}}^1L$ is exact is called an exact Lagrangian.

  If $W$ is compact and the Liouville vector field is strictly outwards pointing at the boundary $\partial W$ and is gradient like for a Morse function $f:W\to \mathbb R$, we say that $W$ is a Weinstein domain. Given a Weinstein domain $W$, one can by a construction called Weinstein handle-attachment attach a Weinstein handle $H_k=D^k\times D^{2n-k}$ to obtain a new Weinstein domain $W\cup H$. All Weinstein domains can be obtained by iterated Weinstein handle-attachments to the unit disc $(D^{2n},\lambda=\frac 12(xdy-ydx))$.

  Given a compatible almost complex structure on $W$, a quadratic complex volume form on $W$ is a trivialization of $(\det_{\mathbb C} T^*W)^{\otimes 2}$. This induces by evaluation a map $\Omega^2:TW^{\times n}\to \C$ which descends to a map $\mathcal LW\to \C/\C^*=S^1$ where $U(n)/O(n)\to \mathcal LW\to W$ is the Lagrangian Grassmannian bundle. Let $L\subset W$ be a Lagrangian, then the map $\mathcal LW\to S^1$ induces a map $F:L\to S^1$. The Maslov index for a loop $\gamma:S^1\to L$ is defined as the degree of the composition $F\circ \gamma:S^1\to S^1$.

  Let $X$ be a Riemannian manifold. Then one can define a Weinstein structure on $D^*X$. The Levi-Civita connection on $D^*X$ induces a canonical splitting $T_{(x,\xi)}D^*X=T_xX\oplus T_xX$ into horizontal and vertical tangent vectors. In this splitting we define the standard compatible almost complex structure $J=\begin{psmallmatrix}0&-\mathrm{id}\\\mathrm{id}&0\end{psmallmatrix}$. We also define the standard quadratic complex volume form on $D^*X$ in the unique way so that $\Omega^2(v_1,\dots, v_n)=1$ for any orthonormal basis of horizontal tangent vectors.

The vanishing of the Maslov class for closed exact Lagrangians in cotangent bundles of closed manifolds is a celebrated result that now has been proven by using several different approaches. The first proof was by Abouzaid  and Kragh using parametrized spectra \cite{MySympfib}. Guillermou proved the result using microlocal sheaves \cite{guillermou2023sheaves}. The most recent proof is by Husin and Kragh using Floer theory \cite{husin2025local}. In \cite{husin2025local} it is furthermore proven that the Maslov class vanishes for closed exact Lagrangians in Weinstein domains obtained from a cotangent bundle of a closed manifold by attaching sub-critical handles over a contractible neighborhood in the base.

In this article we begin by defining cylindrical handles as a generalization of the Weinstein handles defined in \cite{weinstein1991contact}. These cylindrical handles can be attached to a Legendrian link of two spheres in the boundary of a Weinstein domain. We then show that attaching a cylindrical handle is equivalent to first attaching a Weinstein one-handle and then attaching a critical Weinstein handle.

Lastly we prove the following two main theorems. Both theorems show that we cannot remove the assumption that handles are sub-critical in the main result in \cite{husin2025local}.
\begin{theorem}\label{main_result}
  Let $W$ be a Weinstein domain with $\dim W=2n\geq 6$ and with a compatible almost complex structure and a quadratic complex volume form. Then one can to $W$ first attach a Weinstein one-handle $H_1$ and then attach a critical Weinstein handle $H_n$ so that the quadratic complex volume form on $W$ extends to $W\cup H_1\cup H_n$. Furthermore for any extension of the quadratic complex volume form on $W$ to $W\cup H_1\cup H_n$ there exists a closed exact Lagrangian in $W\cup H_1\cup H_n$ with non vanishing Maslov class.
\end{theorem}
The diffeomorphism class of the Lagrangian in Theorem~\ref{main_result} is an $S^{n-1}$ fibration over $S^1$. When $n$ is even this fibration is trivial.
\begin{theorem}\label{main_result2}
  Let $X$ be a closed manifold with $\dim X=n\geq 3$ and assume that there exists an embedded homologically essential sphere $S^{n-1}\subset X$. Then one can to $W=D^*X$ attach a critical Weinstein handle $H$ and find an exact Lagrangian in $W\cup H$ with non vanishing Maslov class with respect to the unique (up to homotopy) extension of the standard quadratic complex volume form on $W$ to $W\cup H$.
\end{theorem}
The Lagrangian in Theorem~\ref{main_result2} has the same diffeomorphism class as $X$ when $n$ is even.
\begin{remark}
  The constructions of the Lagrangians in Theorem~\ref{main_result} and Theorem~\ref{main_result2} works even when $n=2$. However in this case, the statement about the non vanishing Maslov class does not hold.
\end{remark}
In Section~\ref{cylindrical_handles} we provide details on how to construct a cylindrical handle $H\simeq D^*D^1\times D^*S^{n-1}$ and how to attach it to a Legendrian link of two spheres in the boundary of a Weinstein domain. We then show that attaching a cylindrical handle is equivalent to attaching a Weinstein one-handle and a critical Weinstein handle. We argue that the one-handle is attached to a point on each component of the link and that the critical handle is attached to the  connected sum of the link connected by a Legendrian cylinder in the boundary of the one-handle.

In Section~\ref{example} we prove a special case of Theorem~\ref{main_result} where $W=D^{2n}\subset \mathbb R_x^n\times \mathbb R_y^n$ is the $2n$-dimensional unit disc for $n\geq 3$. In this case we attach a cylindrical handle $H$ to the Legendrian link $\{y=0\}\sqcup \{x=0\}\subset \partial D^{2n}$. We then construct two closed exact Lagrangians in $D^{2n}\cup H$ as the two Lagrangian surgeries (introduced by Polterovich in \cite{polterovich1991surgery}) of the singular Lagrangian $\{y=0\}\cup \{x=0\}\subset D^{2n}$, union with a perturbation of the core cylinder $D^1\times S^{n-1}\subset H$. We then compute the Maslov index for a path in these Lagrangians and conclude that they are different and hence the Maslov class of at least one of the Lagrangians is non-zero.

In Section~\ref{local} we prove Theorem~\ref{main_result} in full generality by showing that the construction in Section~\ref{example} can be done locally in any Weinstein domain $W$. In particular, this means that for any cotangent bundle $W=D^*X$ of a closed manifold $X$ with $\dim X\geq 3$, we can attach a cylindrical handle $H$ to a generalized Legendrian Hopf link inside of a contact Darboux ball in $\partial W$. This is equivalent to attaching a Weinstein one-handle and a critical Weinstein handle over a contractible neighborhood in the base. Then for any quadratic complex volume form on $D^*X\cup H$ we can find a closed exact Lagrangian in $D^*X\cup H$ with non vanishing Maslov class. Therefore the assumption that the handles are sub-critical in \cite{husin2025local} is indeed crucial.

In Section~\ref{single} we prove Theorem~\ref{main_result2} by first attaching a critical Weinstein handle $H$ to one of the components of the unit conormal of the homologically essential sphere. We then define an exact Lagrangian in $W\cup H$ as a higher dimensional generalization of the Lagrangian disc mutation (defined in \cite{MR2386535}) of the zero section $X\subset W=D^*X$ with mutation-disc with boundary on $X$ given by the negative Liouville flow applied to the core disc of $H$. This Lagrangian can be described as the zero section $X$ with a neighborhood of the homologically essential sphere replaced by a cylinder that enters the handle. A computation then shows that the Maslov class for that Lagrangian does not vanish. This is another example of why the assumption that the handles are sub-critical in \cite{husin2025local} is required.
\subsection*{Acknowledgements}
The author would like to thank Georgios Dimitroglou Rizell for suggesting the problem and for valuable discussions and comments.
\section{Cylindrical handles}\label{cylindrical_handles}
Let $(W^{2n},\omega=d\lambda)$ be a Weinstein domain. In \cite{weinstein1991contact} Weinstein defined how to perform surgery and attach a handle, topologically $D^k\times D^{2n-k}$, to an isotropic sphere $S^{k-1}$ in the boundary of $W$, in a manner that is compatible with the Weinstein structure and so that the Weinstein structure naturally extends. Note that $k\leq n$ since there are no isotropic submanifolds of dimension higher than $n-1$ in $\partial W$. If $k<n$ the handle is called sub-critical and if $k=n$ the handle is called critical.

In this section we generalize Weinstein's construction and define cylindrical handles, topologically $D^*D^1\times D^*S^{n-1}$. We then explain how to attach a cylindrical handle to a Legendrian link $S^{n-1}\sqcup S^{n-1}\subset \partial W$. Lastly we show that attaching a cylindrical handle is equivalent to first attaching a Weinstein one-handle and then attaching a critical Weinstein handle.

We define a cylindrical handle as
\begin{align*}
  H=D_{r(t)}^*D^1\times D_\delta^*S^{n-1}
\end{align*}
where $D^1=[-1,1]$, $r:D^1\to \mathbb R_{>0}$, $\delta>0$ and $r(t)=\widetilde r>0$ is constant for $t\in [-1,-1/2]\cup [1/2,1]$. The reason we allow the width of the fibers in $D_{r(t)}^*$ to vary with the function $r(t)$ is because we later need enough room for Lagrangians to fit inside the handle. We identify
\begin{align*}
  &D_{r(t)}^*D^1\cong \{(t,y),|t|\leq 1,|y|\leq r(t)\},\\
  &\textstyle D_\delta^*S^{n-1}\cong \{(q_0,\dots,q_{n-1},p_0,\dots,p_{n-1}),\ |q|=1,\ \langle q,p\rangle=0,\ |p|\leq \delta\}.
\end{align*} The symplectic form on $H$ is the standard symplectic form
\begin{align*}
  \omega_H=dt\wedge dy+dq\wedge dp.
\end{align*}
We equip $H$ with a product Liouville form
\begin{align*}
  \lambda_H=\lambda_{D_{r(t)}^*D^1}+\lambda_{D_\delta^*S^{n-1}}
\end{align*} where $\lambda_{D_{r(t)}^*D^1}$ is a Liouville form on $(D_{r(t)}^*D^1, dt\wedge dy)$ and $\lambda_{D_\delta^*S^{n-1}}$ is a Liouville form on $(D_\delta^*S^{n-1},  dq\wedge dp)$. It follows that the Liouville vector field $Z_H$ defined by $\omega_H(Z_H,-)=d\lambda_H$ is also is on product form
\begin{align*}
  Z_H=Z_{D_{r(t)}^*D^1}+Z_{D_\delta^*S^{n-1}}
\end{align*}
where $Z_{D_{r(t)}^*D^1}$ and $Z_{D_\delta^*S^{n-1}}$ are the Liouville vector fields for $\lambda_{D_{r(t)}^*D^1}$ and $\lambda_{D_\delta^*S^{n-1}}$. We require that $Z_{D_{r(t)}^*D^1}$ is gradient like for a Morse function $f_{D_{r(t)}^*D^1}$ and that $Z_{D_\delta^*S^{n-1}}$ is gradient like for a Morse function $f_{D_\delta^*S^{n-1}}$. This implies that $Z_H$ is gradient like for the Morse function $f_H=f_{D_{r(t)}^*D^1}+f_{D_\delta^*S^{n-1}}$. Furthermore we require that
\vspace{3mm}
\begin{enumerate}
\item $Z_{D_{r(t)}^*D^1}$ is inwards pointing at $D_{r(t)}^*D^1|_{t\in \{\pm 1\}}$ and outwards pointing at $S_{r(t)}^*D^1$,
\item $f_{D_{r(t)}^*D^1}$ has a single critical point of index $1$,
\item $Z_{D_\delta^*S^{n-1}}$ is outwards pointing at $S_\delta^*S^{n-1}$,
\item $f_{D_\delta^*S^{n-1}}$ has two critical points, one critical point of index $0$ and one critical point of index $n-1$.
\end{enumerate}
\vspace{3mm}

For instance, we could define
\begin{align*}
  &\lambda_{D_{r(t)}^*D^1}=-ydt +d(\rho(t)y),\\
  &Z_{D_{r(t)}^*D^1}=\rho(t)\partial t+(1-\rho'(t))y\partial y,\\
  &f_{D_{r(t)}^*D^1}=y^2-t^2
\end{align*}
where $\rho:D^1\to \mathbb R$ is a smooth strictly decreasing function such that $\rho(-1)=1$, $\rho(0)=0$, $\rho(1)=-1$ and $|\rho(t)|<\frac{\min r(t)}{\max |r'(t)|}$ for $t\in [-1/2,1/2]$. In this case the critical point of index $1$ to $f_{D_{r(t)}^*D^1}$ is located at the origin. We could also define
\begin{align*}
   &\lambda_{D_\delta^*S^{n-1}}=- pdq+\mu d\left(h(|p|^2)\langle p, \nabla \varphi\rangle \right),\\
  &f_{D_\delta^*S^{n-1}}=|p|^2+\mu h(|p|^2)\varphi
\end{align*}
where $\varphi(q_0,\dots, q_{n-1})=q_0$ is the height function, $\nabla\varphi$ is the gradient with respect to the standard Riemannian metric on $S^{n-1}\subset \mathbb R^n$, $h$ is decreasing and such that $h(t)=1$ for $t$ near $0$ and $h(t)=0$ for $t$ near $1$, and $\mu\in \mathbb R_{>0}$ is so small that the requirements for $Z_{D_\delta^*S^{n-1}}$ and $f_{D_\delta^*S^{n-1}}$ are satisfied. In this case the critical point of index $0$ is located at $(-1,0,\dots,0)\subset S^{n-1}\subset D_\delta^*S^{n-1}$ and the critical point of index $n-1$ is located at $(1,0,\dots,0)\subset S^{n-1}\subset D_\delta^*S^{n-1}$.

By Proposition~4.2 in \cite{weinstein1991contact} there is a neighborhood of the Legendrians
\begin{align*}
  &\{(-1,0)\}\times S^{n-1}\subset D_{r(t)}^*D^1|_{t=-1}\times D_\delta^*S^{n-1}\subset H,\\
  &\{(1,0)\}\times S^{n-1}\subset D_{r(t)}^*D^1|_{t=1}\times D_\delta^*S^{n-1}\subset H
\end{align*} that is isomorphic to the isotropic setup of a standard critical Weinstein handle. Therefore when $\widetilde r$ and $\delta$ are small enough it is possible to attach these ends to each component of a Legendrian link $S^{n-1}\sqcup S^{n-1}\subset \partial W$ and use the Liouville vector field to smoothen the boundary in the same way as a standard handle is attached. The new Weinstein domain obtained by attaching $H$ to $W$ is denoted by $W\cup H$ even though the boundary has been smoothed after taking the union.

By using the standard Riemannian metric on $D^1$ and $S^{n-1}$ we can induce a compatible almost complex structure on the handle. The product of the zero sections $D^1\times S^{n-1}\subset H$ of a cylindrical handle is called the core cylinder.

The Morse function $f_H$ has two critical points, one of index $1$ and one of index $n$. From Lemma~11.13 in \cite{cieliebak2012stein} and the directly following discussion it follows that attaching a cylindrical handle is equivalent to first attaching a Weinstein one-handle $H_1$ to the intersection of the unstable manifold of the critical point of index $1$ with $\partial W$ and then attaching a critical Weinstein handle to the intersection of the unstable manifold of the critical point of index $n$ with $\partial(W\cup H_1)$. By the properties of the Liouville vector field one may show (but we will not need it) that the one-handle will be attached to one point on each component of the link in $\partial W$ and the critical handle will be attached to the connected sum of the link connected by a Legendrian cylinder in $\partial H_1$.

\section{\texorpdfstring{Construction in $D^{2n}\cup H$}{Construction in D2n U H}}\label{example}
In this section we attach a cylindrical handle $H$ to the unit disc $D^{2n}$ for $n\geq 3$ and explicitly construct two closed exact Lagrangians in $D^{2n}\cup H$ using Lagrangian surgery \cite{polterovich1991surgery}, and prove that for any quadratic complex volume form on $D^{2n}\cup H$ one of the Lagrangians must have a non vanishing Maslov class. This proves Theorem~\ref{main_result} in the special case when $W=D^{2n}$.

Let $(D^{2n}\subset \R_x^{n}\times \R_y^n,\omega_{D^{2n}}=d\lambda_{D^{2n}})$ be the unit disc with the Liouville form $\lambda_{D^{2n}}=\frac 12 (xdy-ydx)$. The Liouville vector field $Z_{D^{2n}}=\frac 12(x\partial x+ y\partial y)$ is then outwards pointing at $\partial D^{2n}$. We use the standard compatible (almost) complex structure $J_{D^{2n}}$ and the standard quadratic complex volume form $\Omega_{D^{2n}}^2$ on $D^{2n}$.

Let $0<\epsilon<1/\sqrt{2}$ and let $\psi:D^1\to \R$ be a smooth increasing function defined on $D^1=[-1,1]$ such that $\psi(t)=0$ for $t\leq -\epsilon$, $\psi(0)>0$ and $\psi(t)=t$ for $t\geq \epsilon$. Define exact Lagrangian cylinders $L_1$ and $L_2$ as the image of
\begin{align*}
  &F_1,F_2:D^1\times S^{n-1}\to D^{2n},\\
  &F_1(t,q)= (\psi(t)q, \psi(-t)q),\\
  &F_2(t,q)= (\psi(t)q, -\psi(-t)q)
\end{align*}where $S^{n-1}\subset \R_q^n$ is the unit sphere. The pullback of $\lambda_{D^{2n}}$ to $D^1\times S^{n-1}$ is
\begin{align*}
  &(F_i^*\lambda_{D^{2n}})_{(t,q)}(v)=(\lambda_{D^{2n}})_{F_i(t,q)}\left(\psi(t)v, (-1)^{i-1}\psi(-t) v\right)=0,\quad v\in T_{q}S^{n-1},\\
  &(F_i^*\lambda_{D^{2n}})_{(t,q)}(\partial t)=(\lambda_{D^{2n}})_{F_i(t,q)}\left(\psi'(t)q, (-1)^{i}\psi'(-t)q\right)\\
  &\hspace{2.93cm}=\frac{(-1)^i}2(\psi(t)\psi'(-t)+\psi(-t)\psi'(t))
\end{align*}
so $L_1$ and $L_2$ are indeed exact Lagrangians with primitives $f_i:D^1\times S^{n-1}\to \R$,
\begin{align*}
f_i(t, q)=\frac{(-1)^{i}}2\int_{-1}^t(\psi(\tau)\psi'(-\tau)+\psi(-\tau)\psi'(\tau))d\tau
\end{align*}
satisfying $df_i=F_i^*\lambda_{D^{2n}}$. Note that close to $\partial D^{2n}$ both $L_1$ and $L_2$ coincide with $\{y=0\}\sqcup \{x=0\}$. In particular both $L_1$ and $L_2$ intersect $\partial D^{2n}$ in a generalized Legendrian Hopf link of two $n-1$ dimensional spheres $\Lambda_1=\{y=0\}\cap  \partial D^{2n}$ and $\Lambda_2=\{x=0\}\cap \partial D^{2n}$. The Lagrangians $L_1$ and $L_2$ are known as the two Lagrangian surgeries \cite{polterovich1991surgery} of the singular Lagrangian $\{y=0\}\cup \{x=0\}\subset D^{2n}$.

To the Legendrian link $\Lambda_1\sqcup \Lambda_2\subset\partial D^{2n}$ we attach a cylindrical handle $H=D_{r(t)}^*D^1\times D_{\delta}^* S^{n-1}$, the end $t=-1$ is attached to $\Lambda_1$ by the map
\begin{align*}
  H\supset \{(-1,0)\}\times S^{n-1}&\to \Lambda_1,\\
  ((-1,0),q)&\mapsto (q,0)
\end{align*}and the end $t=1$ is attached to $\Lambda_2$ by the map
\begin{align*}
  H\supset \{(1,0)\}\times S^{n-1}&\to \Lambda_2,\\
  ((1,0),q)&\mapsto (0,q).
\end{align*} On the core cylinder $D^1\times S^{n-1}\subset H$ we define a quadratic complex volume form $\Omega_H^2$ by the requirement that
\begin{align*}
  (\Omega_H)_{(t,q)}^2(\partial t,v_1,\dots,v_{n-1})=e^{i\pi (t+1)(n-2)/2}
\end{align*} for any orthonormal basis $v_1,\dots,v_{n-1}\in T_qS^{n-1}$. This is well defined and continuously extends the standard quadratic complex volume form on $D^{2n}$ and can further, uniquely up to homotopy, be continuously extended to a quadratic complex volume form on all of $D^{2n}\cup H$.

Note that there are other quadratic complex volume forms on $H$ that continuously extend the quadratic complex volume form on $D^{2n}$, but are not homotopic to $\Omega_H^2$. Such quadratic complex volume forms would result in different Maslov indices in the computation below. However the Maslov indices for the paths in the Lagrangians that we will study below would still be different and hence one of these indices would be non-zero.

In $H$ consider the Lagrangians $K_i=dg_i\times S^{n-1}$ where $i\in \{1,2\}$ and  $g_i:D^1\to \R$ is such that $g_i(-1)=-f_i(1, q)$ and $g_i(1)=-f_i(-1,q)=0$ and $g_i'(t)=0$ outside the interval $[-1/2,1/2]$. If these Lagrangians do not fit in the handle, we can make the handle thicker by adjusting $r(t)$ as explained in Section~\ref{cylindrical_handles}. A primitive to the Liouville form $\lambda_H=-ydt +d(\rho(t)y)+\lambda_{D_\delta^*S^{n-1}}$ (defined by the formulas in Section~\ref{cylindrical_handles}) on $K_i$ is $h_i(t,y,q,p)=-g_i(t)+\rho(t)y$.

Since the Lagrangians $K_i$ coincide with the core cylinder $D^1\times S^{n-1}\subset H$ and follow the Liouville flow near the attaching region it follows that $L_i$ and $K_i$ glue together smoothly to closed Lagrangians $L_i\cup K_i\subset D^{2n}\cup H$. Furthermore the primitives $f_i\circ F_i^{-1}:L_i\to \mathbb R$ and the primitives $h_i:K_i\to \mathbb R$ glue together smoothly and hence $L_i\cup K_i\subset D^{2n}\cup H$ are exact.

The Lagrangians $L_i\cup K_i$ have the topology of $S^{n-1}$ fibrations over $S^1$. We continue by computing their Maslov indices when traversing the $S^1$-direction. In $L_1\cup K_1$ consider the path $\gamma_1(t)=(\psi(t)q_0,\psi(-t)q_0)\in L_1$, $t\in D^1$, followed by $\sigma_1(t)=(t,g_1'(t),q_0,0)\in K_1$, $t\in D^1$, here $q_0=(1,0\dots,0)\in S^{n-1}$. The contribution to the Maslov index from $\gamma_1$ with respect to the standard quadratic complex volume form can be computed as the winding number of
\begin{align*}
  \Omega_{D^{2n}}^2(\gamma_1'(t),\partial_{q_2}F_1,\cdots,\partial_{q_n}F_1)
  &=\det\begin{pmatrix}
    \gamma_1'(t)&\partial_{q_2}F_1&\cdots &\partial_{q_n}F_1
  \end{pmatrix}^2\\
                   &=(\psi'(t)-i\psi'(-t))^2(\psi(t)+i\psi(-t))^{2n-2}
\end{align*}
where the columns in the matrix are expressed as complex linear combinations of the basis vectors $\partial x_1,\dots,\partial x_n\in T_{\gamma_1(t)}D^{2n}$. This winding number is $1-\frac n2$. The contribution to the Maslov index from $\sigma_1$ in the handle with respect to $\Omega_H^2$ is computed as the winding number of
\begin{align*}
  \Omega_H^2((1+J g''(t))\partial t,v_1,\dots,v_{n-1})=(1+ig''(t))^2e^{i\pi(t+1)(n-2)/2}
\end{align*}
where $v_1,\dots,v_n$ is any orthonormal basis of $T_{q_0}S^{n-1}$. This winding number is $\frac n2 -1$. In total, the Maslov index for the path $\sigma_1\star \gamma_1\subset L_1\cap K_1$ is $0$.

If we consider the path $\gamma_2(t)=(\psi(t)q_0,-\psi(-t)q_0)\in L_2$, $t\in D^1$, followed by $\sigma_2(t)=(t,g_2'(t),\beta(t),0)\in K_2$ where $t\in D^1$ and $\beta:D^1\to S^{n-1}$ is a path from $q_0$ to $-q_0$, a similar computation as above shows that the Maslov index for $\sigma_2\star \gamma_2\subset L_2\cup K_2$ is $n-2$.

Since the paths $\sigma_1 \star \gamma_1\subset L_1\cup K_1$ and $\sigma_2 \star \gamma_2\subset L_2\cup K_2$ represent the same element in the fundamental group of $D^{2n}\cup H$ and have different Maslov index when $n\geq 3$, no matter how we choose the quadratic complex volume form on $D^{2n}\cup H$ the Maslov index of either $\sigma_1\star \gamma_1$ or $\sigma_2\star \gamma_2$ must be non-zero.
\section{Proof of Theorem 1.1}\label{local}
In this section we prove Theorem~\ref{main_result} by generalizing the result in Section~\ref{example}. We show that one can attach a cylindrical handle $H$ to a generalized Legendrian Hopf link contained in a contact Darboux ball in the boundary of an arbitrary Weinstein domain $W$ with $\dim W\geq 6$. We then show that for any quadratic complex volume form on $W\cup H$ there is a closed exact Lagrangian in $W\cup H$ with non vanishing Maslov class.

Let $(W^{2n},\omega=d\lambda, J, \Omega^2)$, $n\geq 3$ be a Weinstein domain with a compatible almost complex structure $J$ and a quadratic complex volume form $\Omega^2$. Extend the Lagrangians $L_i\subset D^{2n}$ from Section~\ref{example} cylindrically to $\R^{2n}=\R_x^{n}\times \R_y^n$ by
\begin{align*}
  \overline L_i=L_i\cup \{|x|\geq 1,\ y=0\}\cup \{|y|\geq 1,\ x=0\}.
\end{align*}
Choose $z_0\in \partial D^{2n}$ such that the half line $\{tz_0,t\geq 0\}\subset \R^{2n}$ is disjoint from both $\overline L_1$ and $\overline L_2$. For instance we can choose
\begin{align*}
  z_0=\left(\left(\frac 1{\sqrt{2}},0,\dots,0\right),\left(0,\dots, 0,\frac 1{\sqrt{2}}\right)\right).
\end{align*}
Note that such a point $z_0\in \partial D^{2n}$ does not exist when $n=1$. We have that
\begin{align*}
  (\partial D^{2n}\backslash \{z_0\}, \lambda_{\partial D^{2n}})\cong (\R^{2n-1},e^{g(x)}\alpha_{\R^{2n-1}})
\end{align*}
are strictly contactomorphic for some positive function $g:\R^{2n-1}\to \R_{>0}$, here $\lambda_{\partial D^{2n}}=\frac 12(xdy-ydx)$ and $\alpha_{\R^{2n-1}}=dz-ydx$ is the standard contact form on $\R^{2n-1}=\R_x^{n-1}\times \R_y^{n-1}\times \R_z$, to see this we use Proposition~2.1.8 in \cite{geiges2008introduction}. Taking the symplectization of the strict contactomorphism above gives exact symplectomorphisms
\begin{align*}
  (\R^{2n}\backslash \{tz_0,t\geq 0\},\lambda_{\R^{2n}})&\cong (\mathbb R_s\times (\partial D^{2n}\backslash \{z_0\}), e^s\lambda_{\partial D^{2n}})\\
                                        &\cong (\R_s\times \R^{2n-1},e^{s+g(x)}\alpha_{\R^{2n-1}})\\
                                        &\cong (\R_s\times \R^{2n-1},e^{s}\alpha_{\R^{2n-1}})
\end{align*}
where $\lambda_{\R^{2n}}=\frac 12(xdy-ydx)$ and the last map sends $(s,x)$ to $(s+g(x),x)$. When $\epsilon$ from Section~\ref{example} is so small that only the cylindrical parts of $\overline L_1, \overline L_2\subset \R^{2n}\backslash \{tz_0,t\geq 0\}$ are mapped to $\R_{\geq 0}\times \R^{2n-1}$. The images of $\overline L_1$ and $\overline L_2$ result in exact Lagrangian fillings in $(\R_{\leq 0}\times \R^{2n-1}, e^s\alpha_{\R^{2n-1}})$ of a Legendrian link in $\{0\}\times \R^{2n-1}$. By rescaling this filling by a small factor $\beta>0$ in the $x_i$ and $y_i$ directions and by $\beta^2$ in the $z$ direction, we get two exact Lagrangian fillings of a Legendrian link small enough to fit in the symplectization of any choice of contact Darboux ball in the boundary of $W$. We denote these fillings by $\widetilde L_1,\widetilde L_2\subset W$. Now we attach a cylindrical handle $H$ to $\partial \widetilde L_1\sqcup \partial \widetilde L_2\subset \partial W$ and as in Section~\ref{example} extend the Lagrangians $\widetilde L_i$ by $K_i\subset H$ to get closed exact Lagrangians $\widetilde L_i\cup K_i\subset W\cup H$. The computation in Section~\ref{example} shows that for any extension of the quadratic complex volume form on $W$ to $W\cup H$ one of the Lagrangians $\widetilde L_i\cup K_i$ must have a non vanishing Maslov class.
\section{Proof of Theorem 1.2}\label{single}
Assume that $S^{n-1}\subset X$ is a homologically essential sphere in a closed manifold $X$ with $\dim X=n\geq 3$. That $S^{n-1}$ is homologically essential is in these dimensions equivalent to the existence of an embedded closed loop intersecting the sphere transversally exactly once. Let $W=D^*X$ with Liouville form
\begin{align*}
  \lambda_W=-pdq+df
\end{align*}where $f$ is a function with support near the zero section chosen so that $W$ is Weinstein. On $W$ we use the standard compatible almost complex structure $J_W$ and the standard quadratic complex volume form $\Omega_W^2$ with respect to a Riemannian metric on $X$ that we will later specify more details about. In this section we attach a critical Weinstein handle $H$ to one of the components of the unit conormal of the homologically essential sphere. We then construct a closed exact Lagrangian in $W\cup H$ as a higher dimensional generalization of the Lagrangian disc mutation \cite{MR2386535} of the zero section $X\subset W=D^*X$ with mutation-disc with boundary on $X$ given by the negative Liouville flow applied to the core disc of $H$. We then show that the Maslov class for that Lagrangian does not vanish with respect to the unique (up to homotopy) extension of the quadratic complex volume form on $W$ to $W\cup H$.

\begin{figure}[ht]
  \begin{center}
    \begin{tikzpicture}[scale=0.7]
      \draw (-9.5,1.5) -- (9.5,1.5);
      \draw (3,1.5) arc (0:180:3 and 2);
      \draw (-9.5,-1.5) -- (9.5,-1.5);

      \draw (-6,-1.5) -- (-6,1.5);
      \draw (6,-1.5) -- (6,1.5);

      \draw[red] (-9.5,0) -- (-1.5,0);
      \draw[red] (-1.5,0) arc (-90:0:0.5);
      \draw[red] (-1,0.5) -- (-1,1.5);
      \draw[red] (-1,1.5) arc (180:0:1);
      \draw[red] (1,1.5) -- (1,0.5);
      \draw[red] (1,0.5) arc (180:270:0.5);
      \draw[red] (1.5,0) -- (2,0);
      
      \draw[red] (2,0) arc (90:45:1);
      \draw[red] ({2+sqrt(2)/2},{-1+sqrt(2)/2}) -- ({4-sqrt(2)/2},{-3+3*sqrt(2)/2});
      \draw[red] ({4-sqrt(2)/2},{-3+3*sqrt(2)/2}) arc (-135:-45:1);
      \draw[red] ({4+sqrt(2)/2},{-3+3*sqrt(2)/2}) -- ({6-sqrt(2)/2},{-1+sqrt(2)/2});
      \draw[red] (6,0) arc (90:135:1);

      \draw[red] (6,0) -- (9.5,0);

      \node[right] at (1,0.75) {$\gamma_1$};
      \node[above] at (4,{-4+2*sqrt(2)}) {$\gamma_2$};
      \node[above] at (7,0) {$\gamma_3$};
      \node[above] at (-3,0) {$\gamma_3$};
      \node[left] at (-1,0.75){$\gamma_4$};
      \node[above] at (0,2.5) {$\gamma_5$};

      \node[above] at (-5,0) {\color{red}$L$};
      \node[above right] at ({sqrt(2)/2},{1.5+sqrt(2)/2}) {\color{red}$K$};
      
      \node[below left] at (9.5,1.5) {$W=D^*X$};
      \node[below left] at (6,1.5)  {$D^*A$};
      \node[above right] at ({sqrt(2)*3/2},{1.5+sqrt(2)}) {$H$};

      \fill (-6,-1.5) circle (2pt) node[below] {$t=-1$};
      \fill (-3,-1.5) circle (2pt) node[below] {$-1/2$};
      \fill (-1,-1.5) circle (2pt) node[below] {$-\epsilon$};
      \fill (0,-1.5) circle (2pt) node[below] {$0$};
      \fill (1,-1.5) circle (2pt) node[below] {$\epsilon$};
      \fill (3,-1.5) circle (2pt) node[below] {$1/2$};
      \fill (6,-1.5) circle (2pt) node[below] {$1$};
      
      \fill (0,1.5) circle (2pt) node[below] {$\Lambda$};
    \end{tikzpicture}
    \caption{The Lagrangian $L\cup K\subset W\cup H$}\label{f1}
    \end{center}
\end{figure}
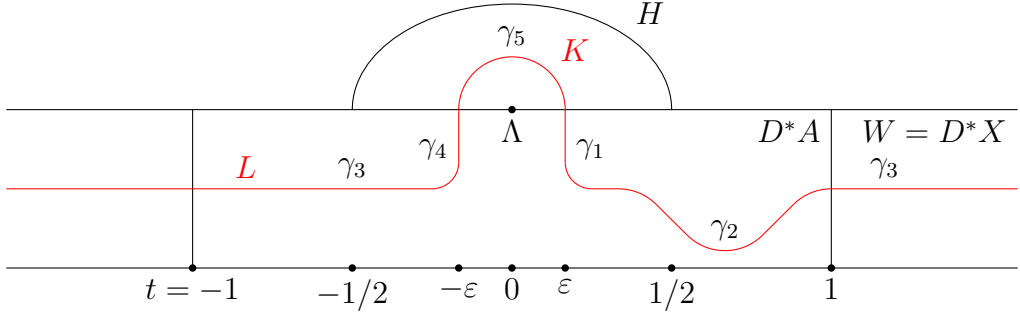

We will now give the details about the constructions which are also schematically illustrated in Figure~\ref{f1}. Take an embedding $A=D_t^1\times S_q^{n-1}\subset X$ that identifies the set $\{0\}\times S^{n-1}\subset A\subset X$ with the homologically essential sphere. Choose a Riemannian metric on $X$ such that the pullback of the metric to $A$ is the product of the standard metrics. Identify
\begin{align*}
  D^1_t\times \mathbb R_y\times T^*S^{n-1}=T^*A\subset T^*X.
\end{align*}
Consider a critical Weinstein handle
\begin{align*}
  &H=D_x^n\times D_{1/2,y}^n,\\
  &\lambda_H=-2ydx-xdy.
\end{align*}
Attach this handle to the Legendrian sphere
\begin{align*}
  \Lambda=\{0\}\times \{1\}\times S^{n-1}\subset S^*A\subset \partial W,
\end{align*}
this is one of the components of the unit conormal of the homologically essential sphere. The point $(0,1,q,0)\in \Lambda$ is attached to $(q,0)\in H$. To make a compatible almost complex structure $J_H$ and a quadratic complex volume form $\Omega_H^2$ on $H$ continuously extend the compatible almost complex structure and quadratic complex volume form on $W$, we define these on the core disc $D^n\times \{0\}\subset H$ as
\begin{align*}
  &J_H(\partial x_i)=\partial y_i,\\
  &\Omega_H^2(\partial x_1,\dots, \partial x_n)=-1.
\end{align*}
We can then uniquely up to homotopy continuously extend the definitions of $J_H$ and $\Omega_H^2$ to get a compatible almost complex structure and a quadratic complex form on all of $W\cup H$.

We now construct a Lagrangian $L\subset W$. Outside of $D^*([-2\epsilon, 1]\times S^{n-1})\subset  D^*A$ we define $L$ to be the zero section. Inside of $D^*([-2\epsilon,1]\times S^{n-1})$ we define $L$ as the union of the images of
\begin{align*}
  &F_1,F_2,F_4:D^1\times S^{n-1}\to D^*A,\\
  &F_1(t,q)=(\epsilon \psi(t)+\epsilon, \psi(-t),q,0),\\
  &F_2(t,q)=((1/2-\epsilon)t+1/2+\epsilon, \theta(t),q,0),\\
  &F_4(t,q)=(-\epsilon \psi(-t)-\epsilon,\psi(t),q,0)
\end{align*}
where $\psi$ is as in Section~\ref{example} and $\theta:D^1\to D^1$ is a function that vanishes near the boundary and such that $\int_{-1}^1 \theta(\tau)d\tau$ have a specific value that will be described later. This specific value of the integral is to ensure that the Lagrangian we end up with in $W\cup H$ is exact. From the equations above it follows that $L$ follows the Liouville flow close to $\partial W$ and that $L$ intersects $\partial W$ at the time $\pm \epsilon$ Reeb flow of $\Lambda$. To make the Lagrangian $L$ exact, we equip the different parts of it with primitives to $\lambda_W$. Outside of $D^*([-2\epsilon,1]\times S^{n-1})$ where $L$ is the zero section, we use the primitive $f$. On the images of $F_1,F_2,F_4$ we use the primitives
\begin{align*}
  &g_1,g_2,g_4:D^1\times S^{n-1}\to \mathbb R,\\
  &g_1(t,q)=\epsilon \int_t^1\psi(-\tau)\psi'(\tau)d\tau+(1/2-\epsilon)\int_{-1}^1\theta(\tau)d\tau+f\circ F_1,\\
  &g_2(t,q)=(1/2-\epsilon)\int_t^1\theta(\tau)d\tau+f\circ F_2,\\
  &g_4(t,q)=-\epsilon\int_{-1}^t\psi(\tau)\psi'(-\tau)d\tau+f\circ F_4
\end{align*}
satisfying $F_i^*\lambda_W=dg_i$.

Consider the Lagrangian $K\subset H$ as the image of $F_5:D^1\times S^{n-1}\to H$
\begin{align*}
  F_5(t,q)=\left (tq,\epsilon\varphi(t) q\right)
\end{align*}
where $\varphi:D^1\to (0,2)$ and $\varphi(t)=t^{-2}$ near the boundary. It follows that the Lagrangian $K$ follows the Liouville flow close to $\partial H$ and that $K$ intersects $\partial H$ at the time $\pm \epsilon$ Reeb flow of the attaching sphere $S^{n-1}\times \{0\}$. Therefore $L$ and $K$ glue together to a smooth Lagrangian $L\cup K\subset W\cup H$. To make the Lagrangian $K$ exact, we equip it with the following primitive to $\lambda_H$
\begin{align*}
  &g_5:D^1\times D^{n-1}\to \mathbb R,\\
  &g_5(t,q)=-\epsilon\int_{-1}^t \left(2\varphi(\tau)+\tau\varphi'(\tau)\right)d\tau-\epsilon\int_{-1}^1\psi(\tau)\psi'(-\tau)d\tau
\end{align*}
satisfying $F_5^*\lambda_H=dg_5$.

If now
\begin{align*}
   \int_{-1}^1 \theta(\tau)d\tau=-\frac {2\epsilon}{1-2\epsilon}\int_{-1}^1(2\varphi(\tau)+\tau\varphi'(\tau))d\tau -\frac{4\epsilon}{1-2\epsilon}\int_{-1}^1\psi(\tau)\psi'(-\tau)d\tau
\end{align*}
the primitives $g_1,g_2,f,g_4,g_5$ glue together smoothly and hence the Lagrangian $L\cup K\subset W\cup H$ is exact.

We will now construct a path $\gamma_5\star\gamma_4\star\gamma_3\star\gamma_2\star\gamma_1\subset L\cup K\subset W\cup H$. Let $q_0=(1,0,\dots,0)\in S^{n-1}$. Let $\gamma_1:D^1\to L\cap D^*A$ be defined by $\gamma_1(t)=F_1(t,q_0)$. Let $\gamma_2:D^1\to L\cap D^*A$ be defined by $\gamma_2(t)=F_2(t,q_0)$. Let $\gamma_3\subset L\cap X$ be a path from $(1,q_0)\in A$ to $(-2\epsilon,-q_0)\in A$, this exists since the sphere $S^{n-1}\subset X$ is homologically essential. Let $\gamma_4:D^1\to L\cap D^*A$ be defined by $\gamma_4(t)=F_4(t,-q_0)$. Finally let $\gamma_5:D^1\to K$ be defined by $\gamma_5(t)=F_5(t,q_0)$.

We will now compute the Maslov index for the loop $\gamma_5\star\gamma_4\star\gamma_3\star\gamma_2\star\gamma_1\subset L\cup K\subset W\cup H$ and show that it is not zero. From Figure~\ref{f1} we see that the contribution to the Maslov index from $\gamma_1$ and $\gamma_4$ are $1/2$ each, this is since the quadratic complex volume form on $D^*A\subset D^*D^1\times D^*S^{n-1}$ is on product form and we only rotate in the $D^*D^1$ factor. The contribution to the Maslov index from $\gamma_2$ is $0$ since one can isotop the Lagrangian through Lagrangians, away from the endpoints of $\gamma_2$, so that it becomes the zero section along the isotoped version of $\gamma_2$. The contribution to the Maslov index from $\gamma_3$ is $0$ since the Lagrangian is the zero section along its path. The contribution to the Maslov index from $\gamma_5$ can be computed as the winding number of
\begin{align*}
                                                                          -\det\begin{pmatrix}
    \gamma_5'(t)&\partial_{q_2}F_5&\cdots &\partial_{q_n}F_5\\
  \end{pmatrix}^2
  &=-\left(1+\epsilon \varphi'(t)i\right)^2\left(t+\epsilon\varphi(t) i\right)^{2n-2}
\end{align*}
which for small $\epsilon$ is close to $1-n$. In the computation above we have used minus the standard quadratic complex volume form on $H$ instead of $\Omega_H^2$. This is good enough since minus the standard quadratic complex volume form and $\Omega_H^2$ coincides on the core disc of $H$, and $K$ lies near that core disc. More precisely the maps these quadratic complex volume forms induce from the Lagrangian Grassmannian bundle on $H$ to $S^1$ differ by an angle less than $\pi/2$ in some neighborhood of the core disc of $H$. Furthermore $K$ lies in that neighborhood when $\epsilon$ is small enough. In total the Maslov index for $\gamma_5\star\gamma_4\star\gamma_3\star\gamma_2\star\gamma_1$ is $2-n$, hence the Maslov class for $L\cup K$ is not zero.

\bibliographystyle{plainurl}
\bibliography{Mybib}

\end{document}